\documentclass[12pt,reqno]{amsart}

\usepackage{amscd}

\headsep=1cm \numberwithin{equation}{section}
 \newtheorem{thm}{Theorem}[section]
 \newtheorem{cor}[thm]{Corollary}
 
 \newtheorem{prop}[thm]{Proposition}
 \theoremstyle{definition}
 
 \theoremstyle{remark}
 \newtheorem{rem}[thm]{Remark}
 \newtheorem*{exam}{Example}
 \numberwithin{equation}{section}
 
\DeclareMathOperator{\depth}{depth}

\DeclareMathOperator{\grade}{grade}
\DeclareMathOperator{\Ass}{Ass}

\newcommand{\fm}{\mathfrak{m}}

\newcommand{\fa}{\frak{a}}

\begin{document}
%-------------------------------------------------------------------------
% editorial commands: to be inserted by the editorial office
%
%\firstpage{1}
%\volume{228}
%\DOI{003-0001}
%\seriesextraline{This is the Concrete Title of this Book\br H.E. R and S.T.C. W, Eds.}
%\firstpage{1}
\title[asymptotic behaviour and Artinian property]
 {asymptotic behaviour and Artinian property of graded local
cohomology modules}
%----------Author 1
\author[S. H. Hassanzadeh]{S. H. Hassanzadeh}

\address{Faculty of Mathematical Sciences and Computer, Tarbiat Moallem University, 599 Taleghani Ave., Tehran 15618, Iran}
\address{Institut de mathematiques, Universite Paris 6, 4, Place
Jussieu, F-75252 Paris Cedex 05, France}

\email{h\_hassanzadeh@tmu.ac.ir}
%----------Author 2
\author[M. Jahangiri]{M. Jahangiri}

\address{Faculty of Basic Science,Tarbiat Modares
University,Tehran,Iran}

\email{m\_jahangiri@modares.ac.ir}
%----------Author 3
\author[H. Zakeri]{H. Zakeri}

\address{Faculty of Mathematical Sciences and Computer, Tarbiat Moallem University, 599 Taleghani Ave., Tehran 15618, Iran}

\email{zakeri@saba.tmu.ac.ir}



\date{30 July 2008}
%%% ----------------------------------------------------------------------

\begin{abstract}
In this paper, considering the difference between the finiteness
dimension and cohomological dimension for a  finitely generated
module, we investigate the asymptotic behavior of grades of
components of graded local cohomology modules with respect to
irrelevant ideal; as long as we study some artinian and tameness
property of such modules.

\end{abstract}

%%% ----------------------------------------------------------------------
\maketitle
%%% ----------------------------------------------------------------------
%\tableofcontents
\section{Introduction}
Throughout the paper, $ R=\bigoplus _{n \geq 0} R_{n}$ is a graded
Noetherian ring where the base ring $ R_{0}$ is a commutative
Noetherian local ring with maximal ideal $ \fm_{0}$ and $R$ is
generated, as an $ R_{0}$-algebra, by finitely many elements of $
R_{1}$. Moreover, we use $ \mathfrak{a}_{0}$ to denote a proper
ideal of $ R_0$ and we set $ R_{+}=\bigoplus _{n>0}R_{n}$, the
irrelevant ideal of $ R$, $\fa=\fa_{0}+R_{+}$, and
$\fm=\fm_{0}+R_{+}$. Also, we use $M$ to denote a non-zero, finitely
generated, graded $R$-module. It is well known (cf \cite[\S12]{BS})
that, for each $i \in \mathbb{N} _ {0}$ (where $\mathbb{N} _ {0}$
denotes the set of all non-negative integers), the $i$-th local
cohomology module $H^{i}_{R_{+}}(M)$ of $M$ with respect to $R_{+}$
inherits natural grading. For each $n\in\mathbb{Z}$ (where
$\mathbb{Z}$ denotes the set of integers), we use the notation
$H^{i}_{R_{+}}(M)_{n}$ to denote the $n$-th graded component of
$H^{i}_{R_{+}}(M)$. Now, according to \cite[15.1.5]{BS}, for each
$i\geq0$, the $ R_{0}$-module $H^{i}_{R_{+}}(M)_{n}$ is finitely
generated  for all $n\in\mathbb{Z}$ and vanishes for all
sufficiently large values of
 $n$. The cohomological dimension of $M$ with respect to $R_{+}$ is
 denoted
 by $cd(R_{+},M)$. Thus $cd(R_{+},M)$ is the largest
non-negative integer $i$ such that $H^{i}_{R_{+}}(M)$ is non-zero.

 The asymptotic behaviors of $H^{i}_{R_{+}}(M)_{n}$ when
$n\to - \infty$ constitute lots of interest (see for example
\cite{B}, \cite{BFL}, \cite{BH}, \cite{BRS} and \cite{M}). As a
basic reference in this topic we recommend \cite{B}. In \cite{BH},
it has been shown that the set
$\Ass_{R_{0}}(H^{f_{R_{+}}(M)}_{R_{+}}(M)_{n})$ is asymptotically
stable, as $n\to - \infty$, where $f_{R_{+}}(M)$, the finiteness
dimension of $M$ with respect to $R_{+}$, is the least non-negative
integer $i$ such that $H^{i}_{R_{+}}(M)$ is not finitely generated.
It therefore follows that
$\grade(\mathfrak{a}_{0},H^{f_{R_{+}}(M)}_{R_{+}}(M)_{n})<\infty$
for all $n\ll0$. Recently, M. Brodmann raised the question whether
the sequence
$(\grade(\mathfrak{a}_{0},H^{f_{R_{+}}(M)}_{R_{+}}(M)_{n}))_{n}$ of
integers is stable, as $n\to - \infty$. One of the main proposes of
this paper is to provide an affirmative answer for this question in
certain cases. In particular, in \ref{2.2}, we show that if
$f_{R_{+}}(M)=cd(R_{+},M)$, then
$\grade(\mathfrak{a}_{0},H^{f_{R_{+}}(M)}_{R_{+}}(M)_{n})=f^{R_{+}}_{\fa}(M)-f_{R
_{+}}(M)$ for all $n\ll0$, where $f^{R_{+}}_{\fa}(M)$, the
$\fa$-finiteness dimension of $M$ relative to $R_{+}$, is the least
non-negative integer $i$ such that $R_{+}$ is not contained in the
radical of the ideal $(0:_{R} H^{i}_{\fa}(M))$.

The concept of tameness is the most fundamental concept related to
the asymptotic behavior of cohomology. A graded $ R$-module $
T=\bigoplus _{n\in \mathbb{Z}} T_{n}$ is said to be tame or
asymptotically gap free (cf \cite[4.1]{BH}) if either $T_n\neq0$ for
all $n\ll0$ or else $T_{n}=0$ for all $n\ll0$. It is known that any
graded Artinian $R$-module is tame \cite[4.2]{B}. Now, the tameness
property of $H^{i}_{\fm_{0}R}(H^{f_{R_{+}}(M)}_{R_{+}}(M))$ provides
some results on the stability of the sequence
$(\depth_{R_{0}}(H^{f_{R_{+}}(M)}_{R_{+}}(M)_{n}))_{n}$, as $n\to
\infty$ (cf corollary \ref{2.6}). In general, the finiteness of
graded local cohomology module has an important role to study the
asymptotic behavior of the $n$-th graded component
$H^{i}_{R_{+}}(M)_{n}$ of $H^{i}_{R_{+}}(M)$, as $n\to - \infty$.
One of parts in finiteness of graded local cohomology is
Artinianess. In the rest of the paper, we study the Artinianess of
the modules $H^{i}_{\fm_{0}R}(H^{j}_{R_{+}}(M))$. Indeed, Theorem
\ref{2.3}, which improves the  results \cite[2.2]{S}
and\cite[4.2]{BRS}, shows that if $j \leq g(M)$, then
$H^{i}_{\fm_{0}R}(H^{j}_{R_{+}}(M))$ is Artinian for $i=0,1$. Also,
Theorem \ref{2.8}, which yields \cite[2.3]{S} as a consequence,
settle, under the assumption $cd(R_{+},M)-f_{R_{+}}(M)=1$, a
necessary and sufficient condition for the $R-$modules
$H^{i}_{\fm_{0}R}(H^{f_{R_{+}}(M)}_{R_{+}}(M))$ to be Artinian.
Finally, \ref{2.10}, establishes some results on the Artinianess of
the $R$-modules $H^{i}_{\fm_{0}R}(H^{cd(R_{+},M)}_{R_{+}}(M))$ in
the case where $dim(R_{0})-3\leq i \leq dim(R_{0})$.

%--------------------------------------------------------------------------------

\section{the results}

We keep the notations and hypotheses introduced in the Introduction.
In addition, for a finitely generated $R_{0}$-module $X$ with
$\fa_{0}X\neq X$, we use the notation $\grade(\fa_{0}, X)$ to denote
the least integer $i$ such that $H^{i}_{\fa_{0}}(X)\neq 0$. In
particular the above mentioned number is denoted by $\depth(X)$
whenever $\fa_{0}=\fm_{0}$.

Let $f_{R_{+}}(M)<\infty$. It is well known (see for example
\cite[5.6]{BH}) that $\Ass_{R_{0}}(H^{f_{R_{+}}(M)}_{R_{+}}(M)_{n})$
is asymptotically stable, as $n\to - \infty$. Therefore, since
$H^{f_{R_{+}}(M)}_{R_{+}}(M)$ is not finitely generated, it follows
that $\fa_{0}H^{f_{R_{+}}(M)}_{R_{+}}(M)_{n}\neq
H^{f_{R_{+}}(M)}_{R_{+}}(M)_{n}$ for all $n\ll0$; hence, for all
$n\ll0$, we have
 $\grade(\mathfrak{a}_{0},H^{f_{R_{+}}(M)}_{R_{+}}(M)_{n})< \infty$.
 Next, we show that if $f_{R_{+}}(M)=cd(R_{+},M)$, then the sequence
$(\grade(\fa_{0},H^{f_{R_{+}}(M)} _{R_{+}}(M)_{n}))_{n}$ is
asymptotically stable, as $n\to - \infty$.

\begin{rem}\label{rem:0}
Using a very slight modification of the proof of \cite[15.1.5]{BS},
one can show that, for all $i\in \mathbb{N}_0$, $H^i_{\fa}(M)=0$ for
all $n\gg 0$. In view of this fact and \cite[2.3]{M} one can see
that
$$f^{R_{+}}_{\fa}(M)=sup\{i:H^j_{\fa}(M)_n=0  \text{~~for all~~}  j<i
 \text{~~and
all~~} n\ll0\}.$$
\end{rem}

With the aid of the above remark and the following proposition, one
can deduce that
$$f^{R_{+}}_{\fa}(M)=inf\{i:H^i_{\fa}(M)_n\neq0 \text{~~for all~~}
  n\ll0\}.$$

%________________________________________________________________________2.2
\begin{prop}\label{2.1}
Let $f^{R_{+}}_{\fa}(M)<\infty$. Then
$H^{f^{R_{+}}_{\fa}(M)}_{\fa}(M)_{n}\neq0$ for all $n\ll0$.
\end{prop}

\begin{proof}
Before beginning the proof, we provide some facts which are needed
in the course of the proof.

Let $y$ be an indeterminate and let $R'_0=R_0[y]_{\fm_0[y]} ,
R'=R'_0 \bigotimes_{R_0}R$ and $M'=R'_0\bigotimes_{R_0}M.$ Then the
natural homomorphism ${R_0}\longrightarrow R'_0$ is faithfully flat
and, by \cite[15.2.2(iv)]{BS}, $H^i_{\fa}(M)_n\bigotimes_{R_0} R'_0
\cong H^i_{\fa R'}(M')_n$ for all $ i\in \mathbb{N}_0$ and all $n\in
\mathbb{Z}$. It therefore follows that $f^{R'_+}_{\fa
R'}(M')=f^{R_+}_{\fa}(M)$. Hence, without loss of generality, we may
assume that $R_0/\fm_0$ is an infinite field. Next we can use the
exact sequences
\[
  H^i_{\fa}(\Gamma_{R_+}(M))\longrightarrow H^i_{\fa}(M)\longrightarrow
 H^i_{\fa}(M/\Gamma_{R_+}(M))\longrightarrow
 H^{i+1}_{\fa}(\Gamma_{R_+}(M))
\]
to see that , for each $i \in \mathbb{N}_0$, $H^i_{\fa}(M)_n \cong
H^i_{\fa}(M/\Gamma_{R_+}(M))_n$ for all but only finitely many $n\in
\mathbb{Z}$. Hence, by replacing $M$ with $M/\Gamma_{R_+}(M)$, we
may assume, in addition, that there exists a homogeneous element $x
\in R_1$ which is a non-zerodevisor on $M$. Now, the exact sequence
$0\longrightarrow M(-1)\stackrel{x}\longrightarrow M \longrightarrow
{M/xM}\longrightarrow 0 $ yields an exact sequence of $R_0$-modules
  \begin{gather}
    H^{i-1}_{\fa}(M)_n\longrightarrow
 H^{i-1}_{\fa}(M/xM)_n\longrightarrow H^i_{\fa}(M)_{n-1}\stackrel{x}\longrightarrow H^{i}_{\fa}(M)_n
  \label{*}\tag{$*$}
  \end{gather}
for all $n\in \mathbb{Z}$ and $i \in\mathbb{N}_0$.Therefore In the
case where $i=f^{R_+}_{\fa}(M)-1$, the exactness of (\ref{*}) in
conjunction with \ref{rem:0} implies that, for all $n \ll0$,
$H^{i-1}_{\fa}(M/xM)_n=0$ ; and hence, by \ref{rem:0}, $
f^{R_+}_{\fa}(M/xM)\geq f^{R_+}_{\fa}(M)-1$.

Now, we prove the assertion by induction on $f=f^{R_+}_{\fa}(M)$.
Let $x \in R_1$ be a non-zero devisor on $M$. Note that there exists
$n_0\in\mathbb{Z}$ such that for all $n \leq n_0$, $
H^{0}_{\fa}(M/xM)_n=0$. Therefore the exact sequence
 (\ref{*}) yields the monomorphism
 \[ \ H^1_{\fa}(M)_{n-1}\stackrel{x}\longrightarrow
   H^{1}_{\fa}(M)_n\]for all $n \leq n_0$.

  If $f=1$,then, by\ref{rem:0}, there exists a descending
sequence of integers, say
  $\{n_i\}_{i\in \mathbb{N}}$, such that $n_i\leq n_0$ and that
 $H^1_{\fa}(M)_{n_i}\neq 0$ for all positive integers $i$. Therefore,
 using the above mentioned monomorphism, one can see that
 $H^1_{\fa}(M)_{n}\neq 0$ for all
 $n \leq n_0$; hence the assertion is true for $f=1$.

 Now, assume inductively that $f \geq 2$ and that the result has
 been proved for smaller values of $f$. Let $M$ be such that
 $f^{R_+}_{\fa}(M)=f$. By the arguments at the beginning of the
 proof, we have $f^{R_+}_{\fa}(M/xM)\geq f-1$. We now consider
 two cases:

 \textbf{Case 1}: $f^{R_+}_{\fa}(M/xM)> f-1$. In this case, in view
 of \ref{rem:0}, there exists $n_0\in\mathbb{Z}$ such that
 $H^{f-1}_{\fa}(M/xM)_{n}= 0$ for all
 $n \leq n_0$. Therefore the exact sequence (\ref{*}) induces the exact
 sequence $ 0\longrightarrow
 H^f_{\fa}(M)_{n-1}\stackrel{x}\longrightarrow
   H^{f}_{\fa}(M)_n$ for all $n \leq n_0$. Now, one can use the
   same arguments as in the case $f=1$ to show that
   $H^{f}_{\fa}(M)_n\neq 0$ for all $n \leq n_0$.

\textbf{Case 2}: $f^{R_+}_{\fa}(M/xM)= f-1$. In this case, by
inductive hypothesis, $H^{f-1}_{\fa}(M/xM)_{n}\neq 0$ for all
 $n \ll 0$.  While, by \ref{rem:0}, $H^{f-1}_{\fa}(M)_{n}= 0$ for all
 $n\ll 0$. Thus,
  using (\ref{*}), we deduce that $H^{f}_{\fa}(M)_{n}\neq 0$ for all
 $n \ll 0$.
The result now follows by induction.
\end{proof}

As an application of the above proposition, we establish, as a
theorem, the following asymptotic behavior of grade.
%_____________________________________________________________________2.3

\begin{thm}\label{2.2} Let $f_\fa^{{R_{+}}}(M)<\infty$ and suppose that
$(f=)f_{R_{+}}(M)=cd(R_{+},M)$. Then, $\grade(\fa_{0},H^{f}
_{R_{+}}(M)_{n})=f^{R_{+}}_{\fa}(M)-f$ for all $n\ll0$.
\end{thm}
\begin{proof}Since $H^{i}_{R_{+}}(M)$ is
finitely generated for all $i<f$, it is easy to see that there
exists $n_{0}\in\mathbb{Z}$ such that, for all $n<n_{0}$,
$H^{i}_{R_{+}}(M)_{n}=0$ for all $i<f$. Therefore, in the
Grothendieck's spectral sequence \cite[11.38]{R}
\[
(E^{p,q}_{2})_{n}=H^{p}_{\fa_{0}R}(H^{q}_{R_{+}}(M))_{n}
\stackrel{p}{\Rightarrow} H^{p+q}_{\fa}(M)_{n}\] we have, for all
$p\in \mathbb{N}_{0}$, $(E^{p,q}_{2})_{n}=0$ for all $n < n_{0}$ and
all $q \neq f$. This in conjunction with \cite[13.1.10]{BS} implies
that $H^{p}_{\fa_{0}}(H^{f}_{R_{+}}(M)_{n})\cong
H^{p+f}_{\fa}(M)_{n}$ for all $n < n_{0}$ and all $p\in\mathbb{
N}_{0}$. Hence, we may use
 \ref{rem:0} and \ref{2.1} to see that, for all $n\ll0$,
$H^{i}_{\fa_{0}}(H^{f}_{R_{+}}(M)_{n})=0$ whenever
$i<f^{R_{+}}_{\fa}(M)-f$ and that
$H^{i}_{\fa_{0}}(H^{f}_{R_{+}}(M)_{n})\neq0$ if
$i=f^{R_{+}}_{\fa}(M)-f$. Thus  $\grade(\fa_{0},H^{f}
_{R_{+}}(M)_{n})=f^{R_{+}}_{\fa}(M)-f$ for all $n\ll0$.
\end{proof}

The asymptotic behavior of $\grade(\fa_{0},H^{j}_{R_{+}}(M)_{n})$,
when $n \to - \infty$ is related to the tameness property of the
graded $R$-modules $H^{i}_{\fa_{0}R}(H^{j}_{R_{+}}(M))$. It is known
that any graded Artinian $R$-module is tame ([1; 4.2]). In the rest
of the paper we look for the cases in which the $R$-modules
$H^{i}_{\fm_{0}R}(H^{j}_{R_{+}}(M))$ are Artinian.

Next, we will use the concept of \emph{cohomological finite length
dimension of $M$} with respect to $R_+$. The concept is defined as
\cite[3.1]{BRS}
\[g( M ) := sup \{ i~: ~\forall j< i, ~~ l_{R_{0}}
(H^{j}_{R_{+}}(M)_{n} ) < \infty ~~\forall n\ll 0~\}.\] The
following theorem  may be viewed as an generalization of
\cite[4.2]{BRS} and \cite[2.2]{S}, where the cases in which $i=0$
and $j=1$ are justified respectively.

%___________________________________________________________________2.4
\begin{thm}\label{2.3} Let $g(M)<\infty$. Then
$H^{i}_{\fm_{0}R}(H^{j}_{R_{+}}(M))$ is an Artinian $R$-module for
$i=0,1$ and $j\leq g(M)$.
\end{thm}
\begin{proof} In view of Kirby's Artinian criterion
for graded modules \cite[Theorem1]{K}, we have to show that, for
$i=0,1$, the
following statements hold.\\
$(i) H^{i}_{\fm_{0}R}(H^{j}_{R_{+}}(M))_{n}$ is an Artinian
$R_{0}$-module for all $n\in\mathbb{Z}$.\\
$(ii) H^{i}_{\fm_{0}R}(H^{j}_{R_{+}}(M))_{n}=0$ for all $n\gg0$.
\\ $(iii) 0:_{H^{i}_{\fm_{0}R}(H^{j}_{R_{+}}(M))_{n}}R_{1}=0$ for
all $n\ll0$.\\   Since, by \cite[13.1.10]{BS} ,
$H^{i}_{\fm_{0}R}(H^{j}_{R_{+}}(M))_{n}\cong
H^{i}_{\fm_{0}}(H^{j}_{R_{+}}(M)_{n})$ for all $i\geq0$ and all
$n\in\mathbb{Z}$, we may use \cite[7.1.3]{BS} and \cite[15.1.5]{BS}
to see that $(i)$ and $(ii)$ hold. So, we only need to prove
$(iii)$. To this end, consider the Grothendieck's spectral sequence
\cite[11.38]{R}
\[(E^{p,q}_{2})_{n}=H^{p}_{\fm_{0}R}(H^{q}_{R_{+}}(M))_{n}
\stackrel{p}{\Rightarrow} H^{p+q}_{\fm}(M)_{n}.\] Using the concept
of $g(M)$, it is easy to see that there exists  $n_{0}\in
\mathbb{Z}$ such that, for all $n<n_{0}$, $(E^{p,q}_{2})_{n}=0$ for
all $q<g(M)$ and all $p\in \mathbb{N}$. Now, the convergence of the
above spectral sequence implies that
$H^{0}_{\fm_{0}R}(H^{j}_{R_{+}}(M))_{n}\cong H^{j}_{\fm}(M)_{n}$ for
all $n<n_{0}$. Therefore, in view of \cite[7.1.3]{BS}, $(iii)$ holds
for $i=0$. While, in the case where $i=1$, the above spectral
sequence yields a monomorphism
\[ H^{1}_{\fm_{0}R}(H^{j}_{R_{+}}(M))_{n} \longrightarrow
  H^{j+1}_{\fm}(M)_{n}
 \] for all $n<n_{0}$. Again, using the fact that
 $H^{j+1}_{\fm}(M)$ is an Artinian graded $R$-module, we have
$0:_{H^{j+1}_{\fm}(M)_{n}}R_{1}
 =0$ for all $n\ll0$; so that,  $0:_{H^{1}_{\fm_{0}R}(H^{j}_
{R_{+}}(M))_{n}}R_{1}
 =0$ for all $n\ll0$, which completes the proof.
\end{proof}

%_________________________________________________________________________________2.5

\begin{rem}\label{2.4}
{\rm The Artinianness of $H^{0}_{\fm_{0}R}(H^{g(M)}_{R_{+}}(M))$ has
already been studied in \cite{BRS}. Also, in \cite{S}, Sazeedeh
shows that $H^{1}_{\fm_{0}R}(H^{1}_{R_{+}}(M))$ is Artinian. As we
mentioned above, these results are consequences of theorem
\ref{2.3}. Moreover, the following example, which has already been
presented in \cite[2.9]{S}, shows that \ref{2.3} is no longer true
for $i=2$.}
\end{rem}

%_________________________________________________________________________________2.6

\begin{exam}\label{2.5}
{\rm Let $k$ be a field and $x,y,t$ be indeterminates. Let
$R_{0}=k[x,y]_{(x,y)}$, $\fm_{0}=(x,y)R_{0}$ and
$R=R_{0}[\fm_{0}t]$, the Ress ring of $\fm_{0}$. Then,
$H^{2}_{\fm_{0}R}(H^{1}_{R_{+}}(R))$ is not Artinian.}
\end{exam}
As an application of \ref{2.3} we have the following corollary.

%_______________________________________________________________________________2.7

\begin{cor}\label{2.6} Let $g(M)<\infty$ and let $j \leq g(M)$. Then
 one of
the following statements hold:\\ $(i)
\depth_{R_{0}}(H^{j}_{R_{+}}(M)_{n})=0$ for all $n\ll0$;\\ $(ii)
\depth_{R_{0}}(H^{j}_{R_{+}}(M)_{n})=1$ for all $n\ll0$;  \\
$(iii) \depth_{R_{0}}(H^{j}_{R_{+}}(M)_{n})\geq2$ for all $n\ll0$.
\end{cor}

\begin{proof}
since , by \ref{2.3}, $H^{i}_{\fm_{0}R}(H^{j}_{R_{+}}(M))$ is tame
whenever $i=0,1$. The result follows immediately.
\end{proof}

%________________________________________________________________________________2.8

\begin{rem}\label{2.7}
{\rm If $f:=f_{R_{+}}(M)=cd(R_{+},M)$, then one can use the same
argument as employed in the proof of \ref{2.2} to see that
$H^{p}_{\fm_{0}}(H^{f}_{R_{+}}(M))_{n}\cong
H^{p+f}_{\fm_{0}+R_{+}}(M)_{n}$ for all $n\ll0$ and all $p\in
\mathbb{N}_{0}$. Hence, using Kirby's Artinian criterion
(\cite[Theorem1]{K}), we deduce that the $R$-module
$H^{p}_{\fm_{0}R}(H^{f}_{R_{+}}(M)$ is Artinian. Note that this
result was obtained in \cite[2.6]{S} under the extra condition that
$R_{+}$ is principal.}
\end{rem}

The next theorem, which is motivated by the above remark, provides a
necessary and sufficient condition for $R$-modules
$H^{i}_{\fm_{0}R}(H^{f_{R_{+}}(M)}_{R_{+}}(M))$ to be Artinian in
the case where $cd(R_{+},M)-f_{R_{+}}(M)=1$. This theorem has been
proved in \cite[2.3]{S} under the further condition that the
arithmetic rank of $R_{+}$ is two.

%____________________________________________________________________________________2.9

\begin{thm}\label{2.8} Let $i\in
\mathbb{N}_{0}$ and let $cd(R_{+},M)=f_{R_{+}}(M)+1$. Then
$H^{i}_{\fm_{0}R}(H^{cd(R_{+},M)}_{R_{+}}(M))$ is an Artinian
$R$-module if and only if
$H^{i+2}_{\fm_{0}R}(H^{f_{R_{+}}(M)}_{R_{+}}(M))$ is an Artinian
$R$-module.
\end{thm}
\begin{proof} Set $f=f_{R_{+}}(M)$ and $c=cd(R_{+},M)$. Consider
the spectral sequence
$(E^{p,q}_{2})_{n}=H^{p}_{\fm_{0}R}(H^{q}_{R_{+}}(M))_{n}\stackrel{p}
{\Rightarrow} H^{p+q}_{\fm}(M)_{n}$. As in the proof of \ref{2.2},
there exists $n_{0}\in \mathbb{Z}$ such that
 $(E^{p,q}_{2})_{n}=0$ for all $n<n_{0}$ and all $q<f$. On
the other hand $(E^{p,q}_{2})_{n}=0$ for all $n\in \mathbb{Z}$ and
all $q>c$. Therefore, the above mentioned spectral sequence induces
an exact sequence of $R_{0}$-modules and $R_{0}$-homomorphisms  \[
0\longrightarrow (E^{i,c}_{\infty})_{n} \longrightarrow
H^{i}_{\fm_{0}R}(H^{c}_{R_{+}}(M))_{n} \stackrel{(d^{i,c}_{ 2})_{n}
}{\longrightarrow} H^{i+2}_{\fm_{0}R}(H^{f}_{R_{+}}(M))_{n}
\longrightarrow (E^{i+2,f}_{\infty})_{n} \longrightarrow 0 \] for
all $i\in \mathbb{N}_{0}$ and all $n<n_{0}$. To prove the assertion,
as in the proof of \ref{2.3}, we employ Kirby's Artinian criterion
for graded modules. It is clear,  by the same argument which is used
in the proof of \ref{2.3}, that the conditions $(i)$ and $(ii)$ of
the criterion are satisfied. So, it is enough for us to verify the
condition $(iii)$. Note that, for each $i,j\in \mathbb{N}_{0}$,
$E^{i,j}_{\infty}$ is a subquotient of the Artinian graded
$R$-module $H^{i+j}_{\fm}(M)$. Therefore if
$H^{i}_{\fm_{0}R}(H^{c}_{R_{+}}(M))$ is Artinian, then, by Kirby's
Artinian criterion, there exists $N\in \mathbb{Z}$ such that
$0:_{(E^{i+2,f}_{\infty})_{n}}R_{1}=0=0:_{im(d^{i,c}_{
2})_{n}}R_{1}$ for all $n<N$. Now, we can use the above exact
sequence to see that
$0:_{H^{i+2}_{\fm_{0}R}(H^{f}_{R_{+}}(M))_{n}}R_{1}=0$ for all
$n<min\{n_{0},N\}$, which, in turn, in conjunction with Kirby's
Artinian criterion implies that
$H^{i+2}_{\fm_{0}R}(H^{f}_{R_{+}}(M))$ is Artinian. The proof of the
reverse implication is similar to the above proof.
\end{proof}

%_______________________________________________________________________________________2.10

\begin{rem}\label{2.9} {\rm Note that, by the same arguments, \ref{2.7}
 and \ref{2.8}
remains true if we replace $f_{R_{+}}(M)$ by $g(M)$.}
\end{rem}
Finally, we get to the edge of the double complex arised from the
spectral sequence
$(E^{p,q}_{2})_{n}=H^{p}_{\fm_{0}R}(H^{q}_{R_{+}}(M))_{n}\stackrel{p}
{\Rightarrow} H^{p+q}_{\fm}(M)_{n}$. The next theorem improves some
of the already known facts for nearer point \cite[2.8]{S} and give
some conditions to investigate the Artinian property for some
farther points.

%________________________________________________________________________________________2.11

\begin{thm}\label{2.10} Set $c=cd(R_{+},M)$ and $d=dim(R_{0})$. Then:\\
$(i) H^{d}_{\fm_{0}R}(H^{c}_{R_{+}}(M))$ and
$H^{d-1}_{\fm_{0}R}(H^{c}_{R_{+}}(M))$ are Artinian $R$-modules.
\\ $(ii)$ If $H^{d-3}_{\fm_{0}R}(H^{c}_{R_{+}}(M))$ and
$H^{d-2}_{\fm_{0}R}(H^{c-1}_{R_{+}}(M))$ are Artinian, then
$H^{d}_{\fm_{0}R}(H^{c-2}_{R_{+}}(M))$ and
$H^{d-1}_{\fm_{0}R}(H^{c-1}_{R_{+}}(M))$ are Artinian.
\\ $(iii)$ If $H^{d-3}_{\fm_{0}R}(H^{c}_{R_{+}}(M))$ is Artinian,
then $H^{d-1}_{\fm_{0}R}(H^{c-1}_{R_{+}}(M))$ is Artinian.
\\ $(iv)$ If $H^{d}_{\fm_{0}R}(H^{c-2}_{R_{+}}(M))$ and
$H^{d-1}_{\fm_{0}R}(H^{c-1}_{R_{+}}(M))$ are Artinian, then
$H^{d-3}_{\fm_{0}R}(H^{c}_{R_{+}}(M))$ is Artinian.
\\ $(v)$ $H^{d-2}_{\fm_{0}R}(H^{c}_{R_{+}}(M))$ is Artinian
if and only if $H^{d}_{\fm_{0}R}(H^{c-1}_{R_{+}}(M))$ is Artinian.
\end{thm}
\begin{proof} Consider the spectral sequence\[
(E^{p,q}_{2})_{n}=H^{p}_{\fm_{0}R}(H^{q}_{R_{+}}(M))_{n}\stackrel{p}
{\Rightarrow} H^{p+q}_{\fm}(M)_{n}.\]  $(i)$. Since $d^{d-1,c}_{
2}=d^{d,c}_{ 2}=d^{d-3,c+1}_{ 2}=d^{d-2,c+1}_{ 2}=0$, we have
$E^{d-1,c}_{\infty}\cong H^{d-1}_{\fm_{0}R}(H^{c}_{R_{+}}(M))$ and
$E^{d,c}_{\infty}\cong H^{d}_{\fm_{0}R}(H^{c}_{R_{+}}(M))$.
Therefore, since, for each $i,j\in \mathbb{N}_{0}$,
$E^{i,j}_{\infty}$ is a subquotient of
 the Artinian graded $R$-module $H^{i+j}_{\fm}(M)$, the assertion
 follows
immediately.
\\ $(ii)\, and \,(iii)$. Consider the exact sequence
\begin{gather}
%\[
0\longrightarrow E^{d-3,c}_{3} \longrightarrow
H^{d-3}_{\fm_{0}R}(H^{c}_{R_{+}}(M))\stackrel{d^{d-3,c}_{
2}}{\longrightarrow} H^{d-1}_{\fm_{0}R}(H^{c-1}_{R_{+}}(M))
\longrightarrow E^{d-1,c-1}_{3} \longrightarrow0.{\tag{\dag}} %\]
\end{gather}
Since $E^{d-1,c-1}_{3}\cong E^{d-1,c-1}_{\infty}$ is Artinian,
$(iii)$ follows immediately from the above exact sequence. Next, as
$E^{d,c-2}_{4}\cong E^{d,c-2}_{\infty}$ and $E^{d-3,c}_{4}\cong
E^{d-3,c}_{\infty}$, we have an exact sequence
\[0\longrightarrow E^{d-3,c}_{\infty}\longrightarrow
E^{d-3,c}_{3}\stackrel{d^{d-3,c}_{ 3}}{\longrightarrow}
E^{d,c-2}_{3}\longrightarrow E^{d,c-2}_{\infty}\longrightarrow 0.\]
Now, we may use the above two exact sequences, together with our
Artinian assumption on $H^{d-3}_{\fm_{0}R}(H^{c}_{R_{+}}(M))$, to
see that $E^{d,c-2}_{3}$ is Artinian. To complete the proof of
$(ii)$, consider the exact sequence
 \[H^{d-2}_{\fm_{0}R}(H^{c-1}_{R_{+}}(M))\stackrel{d^{d-2,c-1}_{
2}}{\longrightarrow} H^{d}_{\fm_{0}R}(H^{c-2}_{R_{+}}(M))
\longrightarrow E^{d,c-2}_{3}\longrightarrow 0\] and use the
Artinian assumption on $H^{d-2}_{\fm_{0}R}(H^{c-1}_{R_{+}}(M))$ to
see that $H^{d}_{\fm_{0}R}(H^{c-2}_{R_{+}}(M))$ is Artinian.
\\ $(iv).$ Note that $E^{d,c-2}_{3}\cong
 \frac{H^{d}_{\fm_{0}R}(H^{c-2}_{R_{+}}(
M))}{im(d^{d-2,c-1}_{ 2})}$. This fact, in conjunction with the
exactness of $ (\dag)$ and our assumptions on
$H^{d}_{\fm_{0}R}(H^{c-2}_{R_{+}}(M))$ and
$H^{d-1}_{\fm_{0}R}(H^{c-1}_{R_{+}}(M))$, implies that
$H^{d-3}_{\fm_{0}R}(H^{c}_{R_{+}}(M))$ is Artinian.
\\ $(v).$ Consider the exact sequence
\[0\longrightarrow E^{d-2,c}_{\infty}\longrightarrow
H^{d-2}_{\fm_{0}R}(H^{c}_{R_{+}}(M))\stackrel{d^{d-2,c}_{
2}}{\longrightarrow}
H^{d}_{\fm_{0}R}(H^{c-1}_{R_{+}}(M))\longrightarrow
E^{d,c-1}_{\infty}\longrightarrow 0
\]and note that $E^{d-2,c}_{\infty}$ and $E^{d,c-1}_{\infty}$ are
 Artinian
$R$-modules. So that the Artinianess of
$H^{d-2}_{\fm_{0}R}(H^{c}_{R_{+}}(M))$ implies the Artinianess of
$H^{d}_{\fm_{0}R}(H^{c-1}_{R_{+}}(M))$ and vice versa.
\end{proof}

{\bf Acknowledgment.} The authors would like to thank Professor M.
Brodmann for raising the question which was discussed in Theorem
\ref{2.2}. \vspace{.3in} \providecommand{\bysame}{\leavevmode\hbox
to3em{\hrulefill}\thinspace}

%-----------------------------------------------------------------------------------------------------------------

% ------------------------------------------------------------------------

%\subsection*{Acknowledgment}
%The authors would like to thank Professor M. Chardin for helpful
%suggestion about Theorem 2.5. Also, the second author would like
%to thank Professor D. Eisenbud for kind consideration to his
%question about unmixed almost
%complete intersection ideals.
% ------------------------------------------------------------------------

\begin{thebibliography}{1}


\bibitem{B} M. Brodmann,\textit{Asymptotic behaviour of cohomology: tameness,
supports and associeted primes}, in: Commutative  Algebra and
Algebraic Geometry, S. Ghorpade, H. Srinivasan, J. Verma, eds.,
Contemp. Math.,\textbf{390}(2005)31-61.

\bibitem{BFL} M. Brodmann and S. Fumasoli and C.S. Lim,\textit{Low codimensional
associeted primes of graded components of local cohomology
modules}, J. Alg.,\textbf{275}(2004)867-882.


\bibitem{BH} M. Brodmann and M. Hellus,\textit{Cohomological pattern of coherent
sheaves over projective schemes}, J. Pure Appl.
Alg.,\textbf{172}(2002) 165-182.


\bibitem{BRS} M. Brodmann, F. Roher and R. Sazeedeh,\textit{Multiplicities of
graded components of local cohomology modules}, J. Pure Appl.
Alg., \textbf{197}(2005)249-278.


\bibitem{BS} M. Brodmann and R.Y. Sharp,\textit{ Local Cohomology: An Algebraic
Introduction with Geometric Applications}, Cambridge Studies in
Advanced Mathematics, 60. Cambridge University Press, Cambridge,
1998.



\bibitem{K} D. Kirby, \textit{Artinian modules and Hilbert polynomials}, Quart. J.
Math. Oxford,  \textbf{24}(2)(1973)47-57.



\bibitem{M} T.Marley, \textit{Finitely graded local cohomology and the depths of algebra}, Proc. Am.
Math. Soc.,  \textbf{123}(1995)3601-3607.



\bibitem{R} J.J. Rotman,\textit{An introduction to homological algebra}, Academic
Press, 1979.


\bibitem{S} R. Sazeedeh, \textit{Finiteness of graded local cohomology modules}, J.
Pure Appl. Alg., \textbf{212}(1)(2008)275-280.




\end{thebibliography}
\end{document}